\documentclass[11pt]{article}
 \usepackage[top=1in,bottom=1in,left=1.5in,right=1.5in]{geometry}
  \usepackage{amsmath,amssymb}
  \usepackage{latexsym}
  \usepackage{graphicx,color}

  \newcommand{\C}{\mathbb{C}}
  \newcommand{\F}{\mathbb{F}}
  \newcommand{\N}{\mathbb{N}}
  
  \newcommand{\R}{\mathbb{R}}

  \newcommand{\U}{\mathbf{U}}
  \newcommand{\uu}{\mathbf{u}}

  \newcommand{\w}{\mathbf{w}}
  
  \newcommand{\x}{\mathbf{x}}
  
  \newcommand{\y}{\mathbf{y}}
  \newcommand{\z}{\mathbf{z}}
  \newcommand{\0}{\mathbf{0}}

  \newcommand{\bF}{\mathbf{F}}
   
  \newcommand{\bH}{\mathbf{H}}

  \newcommand{\cC}{\mathcal{C}}

  \newcommand{\lan}{\langle}
  \newcommand{\ran}{\rangle}
  \newcommand{\an}[1]{\lan#1\ran}
  \def\diag{\mathop{{\rm diag}}\nolimits}
  \newcommand{\hs}{\hspace*{\parindent}}
  \newcommand{\proof}{\hs \textbf{Proof.\ }}

  \newcommand{\trans}{^\top}
  \newcommand{\qed}{\hspace*{\fill} $\Box$\\}

  \newcommand{\dist}{\mathrm{dist}}

  \newcommand{\Sing}{\mathrm{Sing\;}}

  \newcommand{\rB}{\mathrm{B}}
  \newcommand{\rC}{\mathrm{C}}

  \newcommand{\rS}{\mathrm{S}}

  \newcommand{\rank}{\mathrm{rank\;}}

  \newtheorem{theo}{\bfseries \hs Theorem}[section]
  \newtheorem{defn}[theo]{\bfseries \hs Definition}
  \newtheorem{prop}[theo]{\bfseries \hs Proposition}
  \newtheorem{lemma}[theo]{\bfseries \hs Lemma}
  
  \newtheorem{corol}[theo]{\bfseries \hs Corollary}
  
  \newtheorem{example}[theo]{\bfseries \hs Example}
  \newtheorem{rem}[theo]{\bfseries \hs Remark}
  
\newtheorem{notation}{Notation}

  \numberwithin{equation}{section} 

\begin{document}

 \title{Some approximation problems in semi-algebraic geometry}

 \author
 {Shmuel Friedland\footnote{Department of Mathematics, Statistics and Computer Science,
 University of Illinois at Chicago,
 Chicago, Illinois 60607-7045, USA,
 \emph{email}: friedlan@uic.edu. This work was supported by NSF grant DMS-1216393} \ and Ma\l gorzata Stawiska
 \footnote{Mathematical Reviews, 416 Fourth Street, Ann Arbor, MI 48103-4816, USA, \emph{email}: stawiska@umich.edu}}

 \date{January 10, 2016}
 \maketitle

 \begin{abstract}
 In  this paper we  deal with a best approximation of a vector with respect to a closed semi-algebraic set $C$ in the space $\R^n$ 
 endowed with a semi-algebraic norm $\nu$.  Under additional assumptions on $\nu$ we prove semi-algebraicity of the set of points 
of unique approximation and other sets associated 
with the distance to $C$.  For $C$ irreducible algebraic we study the critical point correspondence and  introduce the $\nu$- distance degree, generalizing the notion appearing in \cite{DHOST} for the Euclidean norm. We discuss separately the case of the $\ell^p$ norm ($p>1$).  
 \end{abstract}


 \section{Introduction}\label{sec:intro}

Let $\nu:\R^n\to [0,\infty)$ be a norm on $\R^n$.
 In many applications one needs to approximate a given vector $\x\in\R^n$ by a point $\y$ in a given closed subset $C\subset \R^n$. Usually  the approximation  to $\x$ given by $\y$ is measured by $\nu(\x -\y)$.  Then the distance of $\x$ to $C$ with respect to the norm
 $\nu$ is  $\dist_{\nu}(\x,C):=\min \{\nu(\x-\y): \;\y\in C\}$.   A point $\y^\star\in C$ is called a best $\nu$-approximation of $\x$ 
 if $\nu(\x-\y^{\star})=\dist_{\nu}(\x,C)$.  Let $\|\cdot\|$denote  the Euclidean norm on $\R^n$ and let   $\dist(\x,C)$ denote
 the distance $\dist_{\|\cdot\|}(\x,C)$.  We call a best $\|\cdot\|$-approximation  a best $C$-approximation, or briefly  a best approximation. Of course the Euclidean norm  plays a special role, but it is not the only norm important in practical applications. For example, in compressed sensing one needs to minimize $\|\x\|_1$ subject to linear conditions $A\x=\y$. This motivates us to study approximation problems when the norm $\nu$ is a semi-algebraic function on $\mathbb{R}^n$ and the set $C$ is a closed semi-algebraic subset of $\R^n$. \\

Our paper is organized as follows: In Section \ref{sec:uniqbnuap} we recall the well-known  relation between  the differentiability of $\dist_{\nu}(\cdot,C)$ at a point $\x \in \R^n$ and the uniqueness of $\nu$-approximation of $\x$ (see \cite{Zaj83} for many results of this type, and the references therein). In Section \ref{sec:semalg} we collect some fundamental properties of semi-algebraic sets in $\mathbb{R}^n$ and prove that (under suitable assumptions on the norm $\nu$) the set of all points $\x \in \mathbb{R}^n$ at which the $\nu$-approximation of $\x$ to a (fixed) closed semi-algebraic set  $C$ is not unique is a nowhere dense semi-algebraic set. In subsequent sections we study critical points of the function $\rm{dist } (\x;C)$ given $\x$ and $C$. The interest in critical points is again motivated by applications: In general, the numerical methods for finding best approximation to $\x$ in $C$ aim at finding a local minimum of the function $\nu(\x-\y), \y\in C$, but most of these methods will converge at most to a critical point of $f_{\x,\nu}(\y):=\nu(\x-\y), \y\in C$. In Section \ref{sec:irredvar} we assume  that $C$ is an irreducible algebraic set, but $\nu$ is still a quite arbitrary semi-algebraic norm.  We then study the critical point correspondence $\Sigma(C)$, i.e., the closure of  the semi-algebraic set $\{(\x,\y)\in  \R^n\times C : \y \mbox{ is a critical point of }f_{\x,\nu} \}
\setminus (\R^n\times \Sing C)$. Using some dominating maps associated with stratification of $\Sigma(C)$ we introduce the notion of the $\nu$-distance degree of $C$. The case of $\ell^p$ norms ( $p>1$) is studied in Section \ref{sec:lpnorm}, while Section \ref{sec:algcritptcor} contains  more details about the critical point correspondence. Our results do not require advanced material in semi-algebraic or algebraic geometry, so anyone with basic knowledge of these areas should be able to reconstruct the arguments. At any rate, we include the necessary prerequisites.\\

Some relations with results already present in literature can be observed, but of course there are differences. First of all, the problem of approximating vectors in $\mathbb{R}^n$   by points in a definable closed set $C \subset \mathbb{R}^n$ was studied in \cite{Den11} with respect  to the  Euclidean norm. The author of that paper proved, among other things, that the approximation by points in $C$ is unique outside some nowhere dense definable set. He also obtained an analogous result in the case when $C$ is subanalytic. Both classes considered in \cite{Den11} substantially generalize the class of semi-algebraic sets. However, here we present a straighforward approach to the semi-algebraic case while working with a wide class of semi-algebraic norms. Second, the notions of the critical point correspondence and the Euclidean distance degree of an algebraic variety $C$ were introduced and studied in \cite{DHOST}.  That paper   focused on methods of computation of this degree when $C$ has a parametric representation or when the complexification of $C$ is  an affine cone  in $\mathbb{C}^n$ (which allows one to work in $\mathbb{P}^{n-1}$). Our definition of the distance degree (Definition \ref{defnnudegC})  generalizes that of \cite{DHOST} and works for a quite arbitrary semi-algebraic norm. Finally, in our case the critical point correspondence is a semi-algebraic (rather than algebraic) set and our distance degree agrees with the degree of this semi-algebraic set understood in the sense of \cite{RV02}. We hope that some specialists will take an interest in the natural occurrence of the notion of the degree of a semi-algebraic set in this context. While  we do not present  any explicit computations, we should note that our methods have applications in approximation of matrices and tensors. The relevant  results can be found in \cite{FO12} and \cite{FS13}.

\section{Uniqueness of a best $\nu$-approximation}\label{sec:uniqbnuap}

Let $\nu$ be a norm on $\R^n$.  Let
\begin{equation}\label{defBSnu}
\rB_{\nu}:=\{\x\in \R^n, \nu(\x)\le 1\}, \quad \rS_{\nu}:=\{\x\in\R^n,\nu(\x)=1\},
\end{equation} 
denote respectively the unit ball and the unit sphere with respect to $\nu$.  It is well known that all norms on $\R^n$ are equivalent, i.e.,

\begin{equation}\label{equivnorm}
\kappa_1(\nu)\|\x\|\le \nu(\x)\le \kappa_2(\nu)\|\x\| \textrm{ for all } \x\in\R^n, \textrm{ where } 0<\kappa_1(\nu)\le \kappa_2(\nu).
\end{equation}

Recall that the dual norm $\nu^*$ is defined as $\nu^*(\x)=\max_{\y\in\rS_{\nu}} \y\trans\x$.
Since $\nu$ is a convex function on $\R^n$ it follows that  the hyperplane $\y\trans \z=1, \z\in\R^n$ is a supporting hyperplane of $B_\nu$ at $\x\in\rS_\nu$
if and only if $\y\trans \x=1$ and $\y\in \rS_{\nu^*}$.   
The subdifferential of $\nu$ at $\x\ne \0$ is given by
\begin{equation}\label{defsubdif}\partial\nu(\x):=\{\y\in \rS_{\nu*}, \;\y\trans\x=\nu(\x)\}.
\end{equation} 
Note that $\partial\nu(t\x)=\partial \nu(\x)$ for each $t>0$ and $\x\ne \0$.
In particular, $\nu$ is differentiable at $\x\ne \0$ if and only if the supporting hyperplane of $\rB_\nu$ at $\frac{1}{\nu(\x)}\x$ is unique
\cite{Roc70}.
Assume that $\nu$ is differentiable at $\x\ne \0$.  
Then  the differential of $\nu$ at $\x$, viewed as a row vector in $\R^n$, is the only vector in $\partial\nu(\x)$.  In this case we will denote the differential also by $\partial\nu(\x)$.  
So the directional derivative of $\nu$ at $\x$ in the  direction $\uu\in\R^n$ is given as $ \partial\nu(\x)\uu$.  
We call $\nu$ a differentiable norm if $\nu$ is differentiable at each $\x\ne \0$, i.e., $\nu\in \cC^1(\R^n\setminus \{\0\})$ \cite{Roc70}.
If $\nu$ is a differentiable norm then we denote $\partial\nu(\x)$ by $\nabla\nu(\x)$ for $\x\ne \0$.  

Observe that if $\nu\in \cC^1(\R^n\setminus\{\0\})$ then function $\nu^2\in\rC^1(\R^n)$ satisfied the conditions
\begin{equation}\label{propnusq}
 \nabla \nu^2(\0)=\0  \textrm{ and }\nabla\nu^2(\x)=2\nu(\x)\nabla\nu(\x) \textrm{ for } \x\in
\R^n\setminus\{\0\}.
\end{equation}

Recall that $\nu$ is strictly convex if for each pair of points $\x,\y\in \rS_{\nu}, \x\ne \y$ the point $t\x+(1-t)\y$ lies in the interior of $\rB_{\nu}$ for $t\in (0,1)$,
 i.e. $\nu(t\x+(1-t)\y)<1$
for $t\in (0,1)$.  It is well known that $\nu$ is strictly convex if and only if each supporting hyperplane of $\rB_\nu$ at $\x\in\rS_\nu$ intersects $\rB_\nu$ only in $\x$.
That is, $\nu$ is strictly convex if and only if for each two distinct point $\x_1,\x_2\in\rS_{\nu}$ one has the equality $\partial\nu(\x_1)\cap \partial\nu(\x_2)=\emptyset$.

It is easy to construct a norm in $\R^2$ which is strictly convex and not differentiable.  It is well known that if $\nu$ is differentiable then $\nu^*$ is strictly convex and if $\nu$  is strictly convex then  $\nu^*$ is differentiable.

 Note that the $\ell_p$-norm on $\R^n$, 
$\|(x_1,\ldots,x_n)\trans\|_p=(\sum_{i=1}^n |x_i|^p)^{\frac{1}{p}}$, is differentiable and strictly convex if and only if $p\in (1,\infty)$.
\begin{lemma}\label{realnabnumap}   Assume that $\nu\in \cC^1(\R^n\setminus\{\0\})$.  Then 

\begin{equation}\label{imageRnnu}
\nabla\nu^2(\R^n)=\R^n.
\end{equation}
Suppose furthermore that $\nu$ is a strictly convex norm, i.e. $\nu^*\in \cC^1(\R^n\setminus\{\0\})$.  Then the map $\nabla\nu^2:\R^n\to\R^n$
is one-to-one.
\end{lemma}
\proof
 Assume that $\w\in \R^n\setminus\{\0\}$.  Clearly, $\nabla\nu^2(\w)\ne \nabla\nu^2(\0)=\0$.
Furthermore, there exists a unique supporting hyperplane 
of $B_\nu$ of the form $a(\w)\w$ where $a(\w)>0$.  That is $a(\w)\w\trans \y\le 1$ for all $\y\in B_\nu$ and $a(\w)\w\trans \x=1$ for some $\x\in B_\nu$.
Hence $\nu(\x)=1$ and $\nabla\nu(\x)=a(\w)\w$.  Therefore 

\[\nabla\nu^2(\frac{1}{2a(\w)}\x)=2\nu(\frac{1}{2a(\w)}\x)\nabla\nu(\frac{1}{2a(\w)}\x)=\frac{2}{2a(\w)}\nu(\x)\nabla\nu(\x)=\w.\]  
Thus \eqref{imageRnnu} holds.
Assume that $\nu$ is strictly convex.   Suppose to the contrary that $(\nabla\nu^2)^{-1}(\w)$, for some  $\w\in \R^n\setminus\{\0\}$,  contains at least two distinct points.  These two points are
$\frac{1}{2a(\w)}\x, \frac{1}{2a(\w)}\y$ where $\x,\y\in S_{\nu}$. 
So $\nabla\nu(\x)=\nabla\nu(\y)$ which contradicts the assumption that $\nu$ is strictly convex.\qed 

The following result is known in many variants, dating  back at least to 1938 (\cite{Go38}).  We state and prove here a version best suited to our purposes. For related results, see \cite{Zaj83}, \cite[\S6]{FO12}.

\begin{theo}\label{bestnraprxlemma} Suppose that $\nu$ is strictly convex. Then at each point $\x\not\in C$ where $\text{dist}_\nu(\mathbf{x},C)$ is differentiable, best $\nu$-approximation is unique.
\end{theo}

\proof   The subdifferential of $\nu$ at $\mathbf{a}$ is
\[
\partial \nu (\mathbf{a})=\{\mathbf{x^*} \in (\mathbb{R}^n)^*: \forall \mathbf{y} \ \nu(\mathbf{y}) \geq \nu(\mathbf{a})+\mathbf{x^*}(\mathbf{y}-\mathbf{a})\}.
\]
Let $f(\cdot):= \text{dist}_\nu(\cdot,C)$ Assume that  $f$ is differentiable at $\mathbf{x}$ and let us view the differential of $f$ at  $\mathbf{x}$ as a linear functional on $\mathbb{R}^n$. Let $\mathbf{y^0}$ be a best approximation to $\mathbf{x}$ in $C$. Then $\partial f(\mathbf{x}) \in \partial \nu(\mathbf{x}-\mathbf{y^0})$.\\
If there exists a $\mathbf{z^0} \in C \setminus \{\mathbf{y^0}\}$ such that $\nu(\mathbf{x}-\mathbf{y^0})=\nu(\mathbf{x}-\mathbf{z^0})=f(\mathbf{x})$, then 
$\frac{1}{\nu(\mathbf{x}-\mathbf{y^0})}\partial f(\mathbf{x}) \in \partial \nu(\frac{1}{\nu(\mathbf{x}-\mathbf{y^0})}(\mathbf{x}-\mathbf{y^0}))\cap \partial \nu(\frac{1}{\nu(\mathbf{x}-\mathbf{z^0})}(\mathbf{x}-\mathbf{z^0}))$.\\
But this is a contradiction with strict convexity of $\nu$.   
 \qed

The following  example, due to Sinan G\"unt\"urk (private communication), shows that we cannot drop  the condition that 
$\nu$ is strictly convex  in Theorem \ref{bestnraprxlemma}:
\begin{example} 
 Let $C\subset \R^2$ be the line $\{(t,t), \;t\in\R\}$ and $\nu=\|\cdot\|_1$.  Then for $\x=(x_1,x_2)\trans \in\R^2$
it is straightforward to show that $\dist_{\|\cdot\|_1}(\x,C)=|x_1-x_2|$.  For $\x\in\R^2\setminus C$ a best approximation is an arbitrary point of 
the segment between the two points on $C$: $(x_1,x_1)\trans,(x_2,x_2)\trans$.  That is, for each $\x\in \R^2\setminus C$,  $\dist_{\|\cdot\|_1}(\x,C)$
is differentiable and $\x$ does not have a unique best approximation.

\end{example}

 We will need the following lemma:
\begin{lemma}\label{powernormlem} Let $\nu$ be an arbitrary norm. Let $a>1$ and let $g:\R^n\to \R$ be defined as $g:=\nu^a$.
Then
\begin{equation}\label{gcond}
g(\x)>0 \textrm{ for } \x\ne 0, \quad g(\0)=0, \quad  Dg(\0)=\0, \quad \lim_{\|\x\|\to \infty} g(\x)=\infty,
\end{equation}
where $Dg(\0)$ is the differential of $g$ at $\0$.\\
Assume that $C\subset\R^n$ is a closed set.  Then $\dist_\nu(\x,C)^a$  is differentiable at each $\y\in C$, and its gradient is $\0$. 
\end{lemma} 
\proof Use \eqref{equivnorm} to deduce \eqref{gcond}.  Suppose that $\y\in C$.  Then $\dist_\nu(\y,C)=0$.  Clearly
\[0\le \dist_\nu (\x,C)^a-\dist_\nu(\y,C)^a\le \nu(\x-\y)^a\le (\kappa_2(\nu)\|\x-\y\|)^a.\]
The above inequalities yield directly that $\dist_\nu(\x,C)^a$  is differentiable at each $\y\in C$, and its gradient is $\0$.\qed

Note that in general  the function $\dist(\cdot,C)$ given by a norm $\nu$ may be not differentiable at $\y\in C$. A trivial example is $C=\{\mathbf{0}\} \subset \mathbb{R}^n$.
\begin{corol}\label{coruniqappr} Let the assumptions of Theorem \ref{bestnraprxlemma} hold.  Assume that $a>1$.  
Then at each point $\x$ where $\text{dist}_\nu(\mathbf{x},C)^a$ is differentiable, best $\nu$-approximation is unique.
\end{corol}

\section{Semi-algebraic sets}\label{sec:semalg}

The background material for this section can be found e.g. in \cite{BCR} or \cite{Cos05}.
Recall that a set $S\subset \R^n$ is called semi-algebraic if it is a finite union of sets of the form 
$\{\x \in \mathbb{R}^n: P_i(\x)>0, \ Q(\x)=0, \ i\in \{1,...,\lambda \}\}$, where $P_i, \ Q$ are polynomials on $\mathbb{R}^n$ with real   coefficients. 
The fundamental result about semi-algebraic sets (Tarski-Seidenberg theorem) says that a semi-algebraic set $S\subset\R^n$
can be described by a quantifier-free first order formula (with parameters in $\mathbb{R}$ considered as a real closed field). It follows that the
 projection of a semi-algebraic set is semi-algebraic. The class of semi-algebraic sets is closed under finite unions, finite intersections and complements.\\

A function $f:\R^n \to \R$ is called semi-algebraic if its graph $G(f)=\{(\x,f(\x)): \x\in\R^n\}$ is semi-algebraic.  
The definition and properties of semi-algebraic sets immediately yield the following result:
\begin{lemma}\label{semalgfpa}  Let $f:\R^n\to [0,\infty)$ be semi-algebraic.  Then for each rational $a=\frac{b}{c}$, where $b,c$ are positive integers,
the function $f^a$ is semi-algebraic.
\end{lemma}

The following property is useful:
\begin{prop} (\cite{Shi97}, Theorem I.2.9.13; \cite{Cos05}, Exercise 2.10): Given a semi-algebraic $\mathcal{C}^1$ function $f$ on $\mathbb{R}^n$, the partial derivatives 
$\partial f/\partial x_1,\ldots,\partial f/\partial x_n$ are semi-algebraic.
\end{prop}

From now on we will only consider semi-algebraic norms in $\mathbb{R}^n$. We  call a norm $\nu$ semi-algebraic 
if the function $\nu(\cdot)$ is semi-algebraic.  

\begin{example} The  norm $\|(x_1,\ldots,x_n)\trans\|_a:=
(\sum_{i=1}^n |x_i|^a)^{\frac{1}{a}}, a\ge 1$ is semi-algebraic if $a$ is rational.  Indeed, assume that $a=\frac{b}{c}$, where 
$b\ge c\ge 1$ are coprime integers. Then
\[
G(\|\cdot\|_a)=\{(x_1,\ldots,x_n,t)\trans\! : x_i=\pm y_i^c, y_i\ge 0, i=1,...,n,\  t=s^c,s\ge 0,\sum_{i=1}^n y_i^b-s^b=0\}.
\]
\end{example}

The following result is well known in the case when $\nu$ is the Euclidean norm \cite[\S1.1]{Cos05}. We will sketch the proof in the general case.

\begin{lemma}\label{distsemialg}  Let $C\subset \R^n$ be a nonempty closed semi-algebraic set and let $\nu$ be a semi-algebraic norm. Then the function $f(\cdot):=\dist_\nu(\cdot,C)^a$, where $a=\frac{b}{c}$ with $b,c$  positive integers, is semi-algebraic.
\end{lemma}

\proof  Assume first that $a=b=1$.  Then the graph of $f$ can be written as 
\[
\{(\x,t) \in \mathbb{R}^{n+1}: t\geq 0, \  \forall \y \in C\  t \leq \nu(\x-\y), \ \forall \varepsilon > 0 \ \exists \y^\star\in C: t+\varepsilon > \nu(\x-\y^\star)\}.
\]
 This is a finite intersection of semi-algebraic sets. For example, the set $\{(\x,t) \in \mathbb{R}^{n+1}:  \forall \y \in C \  t \leq \nu(\x-\y)\}$ is 
the complement in $\mathbb{R}^{n+1}$ of the set which is the projection onto $\mathbb{R}^{n+1}$ of the set $B \subset \mathbb{R}^{n+1}\times \mathbb{R}^n$,
\[
B= (\mathbb{R}^{n+1}\times C) \cap u^{-1}(-\infty,0),
\]
where the function $u:\mathbb{R}^{n+1}\times \mathbb{R}^n\mapsto \mathbb{R}, \quad u(\x,t,\y)=\nu(\x-\y)-t$ is semi-algebraic. 
Since preimages of semi-algebraic sets by semi-algebraic maps are semi-algebraic, $B$ is semi-algebraic, and so is its projection.
A similar argument applies to  other sets in the intersection characterizing  the graph of $f$. 
By Lemma \ref{powernormlem} and Lemma \ref{semalgfpa} ,  $\dist_\nu(\cdot,C)^a$ is semi-algebraic if $a=\frac{b}{c}$ and $b,c$ are positive integers.
\qed

For our next theorem, we will need the following results, proved in \cite{Du83} as parts  of Theorem 3.3 (Fact 1 and Fact 2). Recall that a semi-algebraic set is called smooth
if it is an open subset of the set of smooth points of some irreducible algebraic set.  Every semi-algebraic set has  smooth semi-algebraic Whitney stratification, see e.g.  \cite{RV02}.

\begin{prop}\label{prop:strata} Let $X \subset \mathbb{R}^n$ be a semi-algebraic set and let $f: X \mapsto \mathbb{R}^p$ be a semi-algebraic map. Then there is a
smooth Whitney  semi-algebraic stratification $X=\bigcup_{i=1}^L \Delta_i$ 
such that the graph $f\mid \Delta_i$ is a smooth semi-algebraic set for each $i$. 
\end{prop}

\begin{prop}\label{prop:nondiff} Let $X \subset \mathbb{R}^n$ be a smooth semi-algebraic set and let $f: X \mapsto \mathbb{R}^p$ be a map whose graph is 
 a smooth semi-algebraic set. Then the set of points in $X$ where $f$ is not differentiable is contained in a closed semi-algebraic set of dimension
 less than the dimension of $X$.
\end{prop}

Now we can prove an approximation result.  

 \begin{theo}\label{apcon}  Let $C\subset\R^n$ be a closed semi-algebraic set.  Assume that $\nu$ is a semi-algebraic norm such that
$\nu$ and $\nu^*$ are differentiable.
Then the set of all points $\x\in\R^n$ at which the $\nu$-approximation 
 to $\x$ in $C$ is not unique, denoted by $S(C)$, is a nowhere dense semi-algebraic set.
In particular $S(C)$ is contained in some hypersurface $H\subset \R^n$.
 \end{theo}

\proof
Let $f(\x)=\dist_\nu(\x,C)$.  Since $f$ is semi-algebraic, the graph of $G(f)$ is a semi-algebraic set.  Hence 
\begin{equation}\label{defCxGf}
G(f)\times C:=\{(\x\trans,t,\y\trans)\trans \in\R^{2n+1},\; \x\in\R^n, \ t=\dist_\nu(\x,C), \ \y\in C\}
\end{equation}
is semi-algebraic.  Let 
\begin{eqnarray}\notag
T(f):=&&\{{(\x\trans,t,\y\trans,\x\trans,t,\z\trans)}\trans \in\R^{2(2n+1)},
 (\x\trans,t,\y\trans)\trans,(\x\trans,t,\z\trans)\trans\in G(f)\times C,\\
&&\nu(\x-\y)=\nu(\x-\z)=t,\;\|\y-\z\|^2>0\}.\label{defSing1}
\end{eqnarray}
Clearly, $T(f)$ is semi-algebraic. It is straightforward to see that $S(C)$ is the projection of $T(f)$ on the first $n$ coordinates, so  $S(C)$ is semi-algebraic.  
Theorem \ref{bestnraprxlemma}  along with Propositions \ref{prop:strata} and \ref{prop:nondiff} yields that $S(C)$  does not contain an open set.   Therefore $S(C)$ is contained
in a finite union of hypersurfaces, which is a hypersurface $H$.  \qed

The proof of Theorem \ref{apcon} yields 
\begin{corol}\label{defOmegnuC}
Let the assumption of Theorem \ref{apcon} hold.  Define
\[\Omega_\nu(C):=\{(\x,\z), \;\x\in\R^n\setminus S(C), \z\in C, \dist_\nu(\x,C)=\nu(\x-\z)\}.\]
Then $\Omega_\nu(C)$ is a semi-algebraic set of dimension $n$.
\end{corol}

\begin{rem} A similar result for the Euclidean norm and a definable or subanalytic set was proved in \cite{Den11}. The proof uses strict convexity of the norm in an essential way.  We want to highlight the fact that our approach implies easily the fact that $S(C)$ is semi-algebraic. In a general context (of metric geometry, not necessarily semi-algebraic) there are many ways of proving that $S(C)$ is nowhere dense, see e.g. \cite{Zaj83}. Also, in \cite{Fre97} it is proved directly that the dimension of $S(C) \subset \mathbb{R}^n$ is at most $n-1$ (the argument works for a few standard notions of dimension, e.g. the Hausdorff dimension). We thank the referee for the latter reference. From the definition of dimension of a semi-algebraic set  it follows in a straightforward way that a semi-algebraic set $S$ has an empty interior if and only if ${\rm dim }\; S<n$.
\end{rem}

\section{The case of an irreducible variety}\label{sec:irredvar}

We first recall some basic facts about varieties and polynomial and rational maps used in this paper, see for example \cite{Cos05, Mi68, Mu95, Na66}.
We will work only with real and complex algebraic varieties, so for the purpose of stating general results we let $\F$ denote either the field of real numbers $\R$ or complex numbers $\C$.
Let $\F[\F^n]$ be the ring of polynomials in $n$ variables $\x=(x_1,\ldots,x_n)\trans\in\F^n$.
For $p\in \F[\F^n]$ denote by $Z(p)\subset \F^n$ the zero set of $p$.
$V\subset \F^n$ is called a variety or an algebraic set if there exists a finite number of polynomials $p_1,\ldots,p_m\in \F[\F^n]$
so that $V=\bigcap_{i=1}^m Z(p_i)$.  Let $V\subset \F^n$ be a variety and  let $I_V\subset \F[\F^n]$ be the ideal of all polynomials
that vanish on $V$.  Then by Hilbert's theorem $I_V$ is  finitely generated.  We assume here that $I_V$ is generated by $p_1,\ldots,p_m$,
i.e. $I_V=\an{p_1,\ldots,p_m}$. 
A variety $V$ is called  irreducible  if $V$ is not a union of two proper subvarieties of $V$.  
Recall that $V$ is irreducible if and only if $I_V$ is prime.
It is well known every variety $V$ decomposes uniquely as a finite union of distinct irreducible varieties.

For two given varieties $X,Y\subset \F^n$, the set $X\setminus Y$ is called a quasi-variety.  It is well known that the closure of a quasi-variety in 
a standard or Zariski topology is a variety. 
For $S\subset \F^n$ we denote by Closure$(S)$ the closure of $S$ in the standard topology in $\F^n$.

Assume that $V\subset \F^n$ is a variety defined as above.  Denote by $D(V)(\x)=([\frac{\partial p_i}{\partial x_j}(\x)]_{i=j=1}^{m,n})\trans\in \F^{n\times m}$, 
\  $\x=(x_1,\ldots,x_n)\trans \in \F^n$, the Jacobian matrix corresponding to $V$.  Then the dimension of $V$, denoted by $\dim V$, is the minimal
possible nullity of $D(V)(\y)$ for $\y\in V$. (Recall that for an $n\times m$ matrix $A$, the nullity of $A$ is the dimension of the null space of $A$, i.e., of  $\{\x: A\x=\0\}$.) 
Assume that $V\subset \F^n$ is an irreducible variety of dimension $d, 1\le d<n$.  A point $\y\in V$ is called smooth if
$\rank D(V)(\y)=n-d$.  Otherwise $\y\in V$ is called singular.   The set of singular points of $V$ is denoted by 
$\Sing V$. Note that $\Sing V$ is the set of all points of $V$ where all $n-d$ minors of $D(V)(\y)$ are zero, 
hence $\Sing V$ is a strict subvariety  of $V$.  The quasi-variety
$W:=V\setminus \Sing V$ is a nonempty manifold of dimension $d$.  For $\F=\C$ it is connected.  For $\F=\R$ it consists of a finite number of connected components.
For each $\y\in V\setminus \Sing V$ we denote by $\U(\y)\in\F^n$ the $n-d$-dimensional subspace spanned by the columns of $D(V)(\y)$. 

More generally, assume that the quasi-variety $W\subset \C^n$ is a complex connected manifold of dimension $d$.  Then its closure
is an irreducible variety.  

Suppose that $C=\bigcap_{i=1}^m Z(p_i)\subset \R^n$ is an irreducible variety of real dimension $d$, where 
$I_C=\an{p_1,\ldots,p_m}$.  Then $C_\C$ is the zero locus of $p_1(\z)=\ldots=p_m(\z)=0$
in $\C^n$.  It is well known that $C_\C$ is a complex irreducible variety of complex dimension $d$, see e.g. \cite{RV02}.
We denote $D(C_\C)(\x)$ by $D(C)(\x)$ when no ambiguity would arise.

Let $X,Y$ be affine irreducible algebraic varieties over $\C$.  A map  $f: X \mapsto Y$ is a regular map if it is polynomial in affine coordinates on $X,Y$.
A map $f: X \dashrightarrow Y$ is called rational if there exists a Zariski open set $X'\subset X$ such that the restriction $f:X'\mapsto Y$ is given
as a well defined rational map in affine coordinates on $X,Y$. 
Such a map $f$ is called dominant if $f(X')$ is Zariski dense in $Y$.

Let $X,Y$ be irreducible affine varieties of the same dimension and let $f: X \mapsto Y$ be a regular dominant map. Then the degree of $f$, denoted by $\deg f$, 
is defined as the (necessarily finite) degree of the field extension $[\mathbb{C}(X):f^*(\mathbb{C}(Y))]$.
Furthermore, $\deg f$ is the cardinality of the set $f^{-1}(\y)$ for a generic $\y\in Y$.

The following result, introducing the notion of a critical point, is well known, see for example Lemma 2.7 in \cite{Mi68} for a polynomial real-valued function $g$, and we leave its proof to the reader:
\begin{lemma}\label{charcritpoits}  Let $C\subset \R^n$ be an irreducible variety.  
Assume that $g\in \cC^1(\R^n)$.  Then $\y\in C$ is a critical pointof $g\mid C$  if one of the following equivalent conditions holds:
\begin{enumerate} 
\item Either $\y\in \Sing C$ or $\y\in C\setminus \Sing C$ and $\nabla g(\y)\mid T_{\y}C \equiv 0$, where $T_{\y}C$ is the tangent space to $C$ at $\y$. 
\item Either $\y\in \Sing C$ or $\y\in C\setminus \Sing C$ and $\nabla g(\y)\in \U(\y)$, where $\U(\y)$ is the column space of $D(C)(\y)$, i.e., the normal space to $C$ at $\y$.
\item  Either $\y\in \Sing C$ or $\y\in C\setminus \Sing C$ and $\rank [D(C)(\y)\;\nabla g\trans]=n-d$.
\item $\rank  [D(C)(\y)\;\nabla g\trans]\le n-d$.
\item $\y\in C_\C$ is in the zero set of all $n-d+1$ minors of $ [D(C)(\y)\;\nabla g\trans]$.
\end{enumerate}
\end{lemma}

We will now study the properties of the set of critical points of $g$.  As the singular points of $C$ are always critical points of $g$ by the definition, it is natural 
to consider only the smooth points of $C$ which are critical points $g$ \cite{Mi68}.

\begin{notation}\label{edfnprojpi} Let $m>n$ be two positive integers and $\F$ be the field of real $\R$ or complex numbers $\C$.   For the remaining part, we  let $\pi:\F^m \to \F^n$ denote the projection on the $n$ first coordinates of vectors in $\F^m$.
\end{notation}

\begin{lemma}\label{theo:CritSemialg} Let $g\in\cC^1(\R^n)$ be semi-algebraic.  For each $\x\in\R^n$ denote by $g_\x:\R^n\to \R$
the function $g_\x(\y)=g(\x-\y)$.  Assume that $C\subset\R^n$ is an irreducible variety.
Then the sets: 
\begin{eqnarray*}
&&\Sigma_{g,0}(C) =\{(\x,\y)\in  \R^n\times C : \y \mbox{ is a critical point of } g_{\x}\mid C\},\\
&&\Sigma_{g,1}(C):=\Sigma_{g,0}(C)\setminus (\R^n\times \Sing C),\\
&&\Sigma_g(C):=\textrm{Closure } (\Sigma_{g,1}(C))
\end{eqnarray*}
are semi-algebraic. 
Then $\pi(\Sigma_{g,1}(C))$ is semi-algebraic.
\end{lemma} 
\proof
$\Sigma_{g,0}(C)$ is  the complement   of the projection (onto the product of the first two factors) of the semi-algebraic set
\[
B=\{(\x,\y,\z) \in  \R^n\times \R^n \times \R^n: \y \in C\setminus \Sing C, \ \z \in T_{\y}C, \ \nabla g_{\x}(\y)(\z) \neq 0 \}.
\]
Hence $\Sigma_{g,0}(C)$ is semi-algebraic.  Clearly, $\R^n\times\Sing C$ is an algebraic set.  
Hence $\Sigma_{g,1}(C)$ is semi-algebraic, and its closure $\Sigma_g(C)$ is semi-algebraic.\qed

Suppose that $g$ is a polynomial.  This will be the case when e.g. 
$g=\|x\|_q^q$ with $q \geq 2$ even. Then we can define the sets $\Sigma_{g,0}(C_\C), \Sigma_{g,1}(C_\C), \Sigma_g(C_\C)$ over $\C$ as in Lemma \ref{theo:CritSemialg}, where $\R^n,C,\Sing C$ are replaced by $\C^n,C_{\C},\Sing C_{\C}$  respectively. The arguments of the proof of
Lemma \ref{charcritpoits} yield that $\Sigma_{g,0}(C_\C)$ is a complex subvariety.  Hence $\Sigma_{g,1}(C_\C)$ is a quasi-variety and $\Sigma_g(C_\C)$
is a variety. See  \cite{FS13} for more details.
\begin{lemma}\label{interseclem}  Suppose that $\nu$ is a semi-algebraic norm such that  $\nu$ and $\nu^*$ are differentiable. 
Assume that  $g=\nu^{a}$, where $a=\frac{b}{c}$ and $b>c\ge 1$ are integers.
Let $C\subset \R^n$ be an irreducible variety and assume that $\Omega_\nu(C)$ is defined as in Corollary \ref{defOmegnuC}.  
Then $\Omega_\nu(C)\subseteq \Sigma_{g,0}(C)$.
Furthermore, $\Omega_\nu(C)\setminus (\R^n\times \Sing C)=\Omega_\nu(C)\cap \Sigma_{g,1}(C)$,  and  
$\Omega_\nu(C)\setminus (\R^n\times \Sing C),\pi(\Omega_\nu(C)\setminus (\R^n\times \Sing C))$ 
are semi-algebraic sets of dimension $n$.
\end{lemma}
\proof   Suppose first that $\x\in C$.  Then $\x$ is a unique best $\nu$-approximation of $\x$ and $(\x,\x)\in\Omega_\nu(C)$. 
As $\nabla \dist(\cdot,C)^a$ is $\0$ at $\x$, we deduce that $(\x,\x)\in \Sigma_{g,0}(C)$.
Let $\x\in\R^n\setminus (S(C)\cup C)$.   
So there is a unique $\z$ such that $\dist_\nu(\x,C)=\nu(\x-\z)>0$, i.e.,  $(\x,\z)\in \Omega_\nu(C)$. 
If $\z\in\Sing C$ then $(\x,\z)\in \R^n\times \Sing C\subset  
\Sigma_{g,0}(C)$.  Assume now that $\z\in C\setminus \Sing C$.  Since $\dist_\nu(\x,C)=\nu(\x-\z)>0$ it follows that $\nabla g_\x(\z)=0$.  Hence
$(\x,\z)\in \Sigma_{g,0}(C)$ and $\Omega_\nu(C)\subseteq \Sigma_{g,0}(C)$.  Furthermore, 
$\Omega_\nu(C)\setminus (\R^n\times \Sing C)=\Omega_\nu(C)\cap \Sigma_{g,1}(C)$.
As $\Omega_\nu(C)$ is semi-algebraic, it follows that  $\Omega_\nu(C)\setminus (\R^n\times \Sing C)$ and $\pi(\Omega_\nu(C)\setminus (\R^n\times \Sing C))$
are semi-algebraic.

Let $\y\in C\setminus\Sing C$.  Then $\min_{\w\in\Sing C}\nu(\y-\w)=\phi(\y)>0$.  Let  $O(\y):=\{\x\in \R^n,\;\nu(\x-\y)<\frac{\phi(\y)}{2}\}$.
Note that $O(\y)$ is an open semi-algebraic set, hence its dimension equals $n$.
Assume that $\x\in O(\y)$.
Then $\dist_\nu(\x,C)=\nu(\x-\z)\le \nu(\x-\y)< \frac{\phi(\y)}{2}$.  Suppose that $\w\in \Sing C$.  Then 
\[\nu(\x-\w)=\nu(\y-\w+\x-\y)\ge \nu(\y-\w)-\nu(\x-\y)>\phi(\y)-\frac{\phi(\y)}{2}=\frac{\phi(\y)}{2}.\] 
This shows that $\z\not\in \Sing C$.  Thus for each $\x\in O(\y)\setminus S(C)$ we have that $(\x,\z)\in \Omega_\nu(C)\setminus (\R^n\times \Sing C)$.
Hence $\dim \Omega_\nu(C)\setminus (\R^n\times \Sing C)=\dim \Omega_\nu(C)=n$.  Clearly, 
 $\pi(\Omega_\nu(C)\setminus (\R^n\times \Sing C))\supset O(\y)\setminus S(C)$.  Hence  $\pi(\Omega_\nu(C)\setminus (\R^n\times \Sing C))$ is a semi-algebraic set of
dimension $n$.\qed

Assume that $\nu(\cdot)=\|\cdot\|$ and $g(\cdot)=\|\cdot\|^2=\sum_{i=1}^n x_i^2$.  Then $\Sigma_g(C_\C)$ is an irreducible variety of dimension $n$ \cite{DHOST, FS13}.
(A short justification that $\Sigma_g(C_\C)$ is irreducible is given in the proof of Lemma \ref{sig2malgvar}.) 
$\Sigma_g(C_\C)$ is called the critical point correspondence  in \cite{DHOST}.   Consider the projection $\pi:\Sigma_g(C_\C)\to \C^n$.  
Clearly, $\pi$ is a polynomial map.  It is a dominating map \cite{DHOST, FS13}.
This is a simple consequence of Lemma \ref{interseclem}.  Indeed, since  $\pi(\Omega_\nu(C)\cap \Sigma_{g,1}(C))$ is a semi-algebraic set of real dimension $2n$
it follows that the algebraic set Closure$(\pi(\Sigma_g(C_\C))$ must be $\C^n$. As the complex dimension $\dim \Sigma_g(C_\C)=n$, it follows that $\pi$ is a dominating map.\\

Let $\delta:=\delta_{\|\cdot\|}(C)$ be the degree of  $\pi:\Sigma_g(C_\C)\to \C^n$.  That is, for a generic point $\x\in\C^n$ the set $\pi^{-1}(\x)\subset \Sigma_g(C_\C)$
has $\delta$ distinct points:
\[\pi^{-1}(\x)=\{(\x,\z_1),\ldots,(\x,\z_{\delta})\}.\]
The number $\delta$ is called the Euclidean distance degree of $C$ in \cite{DHOST}.   It gives an upper bound for the number of critical smooth points for the function
$g_\x\mid C$ for a generic $\x\in \R^n$.
\begin{theo}\label{theodefdegnuC}  Let the assumptions of Lemma \ref{interseclem} hold.
Assume that $\Omega_\nu(C)\setminus (\R^n\times \Sing C)=\bigcup_{i=1}^N \Theta_i$ is a 
  smooth semi-algebraic Whitney stratification.  Suppose that each $\Theta_i$ is an open semi-algebraic subset of smooth points of an 
irreducible variety $V_i\subset\R^n\times\R^n$ 
of dimension $n$ for $i=1,\ldots,M$. Furthermore, for $i>M$ each $\Theta_i$ has dimension less than $n$.  Denote by $V_{i,\C}\subset \C^n\times \C^n$ 
the complexification of $V_i$.  Then the projection $\pi: V_{i,\C}\to \C^n$ is dominating for $i=1,\ldots,M$.   
\end{theo}
\proof  Clearly $\Omega_\nu(C)$ is a univalent graph and
\[\Theta_i=\{(\x,\z), \x\in \pi(\Theta_i), \z\in C\setminus\Sing C \textrm{ is the unique best approximation of } \x\}.\]
Hence $\dim \Theta_i=\dim \pi(\Theta_i)$.   Assume that $i\in\{1,\ldots,M\}$.  Then $n=\dim V_{i,\C}\ge \dim \pi(\textrm{Closure} (\pi(V_{i,\C}))$.
As $\pi(\Theta_i)\subset \pi(V_{i,\C})$ and $\dim \pi(\Theta_i)=n$ it follows that $\pi: V_{i,\C}\to \C^n$ is dominating.\qed

We now give a definition of the $\nu$-distance degree of $C$:
\begin{defn}\label{defnnudegC}  Let the assumptions of Theorem \ref{theodefdegnuC} hold.  
Let $\delta_i$ be the degree of the dominating map $\pi:V_{i,\C}\to \C^n$ for $i=1,\ldots,M$.  
Assume that $V_1,\ldots,V_k$, where $k\le M$,
are $k$ distinct irreducible varieties, while  $V_i\in\{V_1,\ldots,V_k\}$ for each $i>k$.  
Then $\delta_\nu(C):=\sum_{i=1}^k \delta_i$ is called the $\nu$-distance degree of $C$.
\end{defn}

\begin{rem}
The union $\bigcup_i V_{i,\mathbb{C}}$ is the smallest complex algebraic variety containing $\Omega_\nu(C)\setminus (\R^n\times \Sing C)$. Therefore the $\nu$-distance degree agrees with the degree of the semialgebraic set $\Omega_\nu(C)\setminus (\R^n\times \Sing C)$ as introduced in \cite{RV02}. 
Note also that for $\nu(\cdot)=\|\cdot\|$ our definition of $\delta_{\|\cdot\|}(C)$ coincides with the Euclidean distance degree of $C$ in \cite{DHOST}. There are many tools of algebraic geometry (including computational ones) that can be used to find the degree of the dominating map $\pi: V_{i,\C}\to \C^n$
\cite{DHOST, FS13, ChiI, ChiII, ChiIII}.
\end{rem}

In the last section we interpret the degree of $\pi:V_{i,\C}\to \C^n$ as a degree of critical point correspondence of an algebraic function 
induced by one stratum of dimension $n$ given by Proposition \ref{prop:strata} for the graph of the norm $\nu$ considered as a subset of $\mathbb{R}^n\times \mathbb{R}$.
\section{Critical point correspondence for  $\ell_p$ norms}\label{sec:lpnorm}

Let $m\in \N$ and let $\bF_{2m-1}:\C^n\to\C^n$ denote the map $(z_1,\ldots,z_{n})\mapsto (z_1^{2m-1},\ldots,z_n^{2m-1})$.
Note that $\bF_{2m-1}$ is a proper map of degree  $(2m-1)^n$.  For each irreducible variety $W\subset \C^n$ of dimension $t$, $\bF^{-1}(W)$ is a union of at most
$(2m-1)^n$ irreducible varieties, each of dimension $t$.  Let $X\subset \C^n$ be an irreducible variety of dimension $d$.   Then
$\bF_{2l+1}(X)$ is also an irreducible variety of dimension $d$.

Assume that $p=\frac{2m}{2l+1}$, where $m>l\ge 0$ are integers.  Then 
\[
\|\x\|_p=\big(\sum_{j=1}^n x_i^{\frac{2m}{2l+1}}\big)^{\frac{2l+1}{2m}}, \ \x \in \mathbb{R}^n.
\]
Note that $\|\x\|_p$ is a $C^\infty$-smooth  function on $\R^n\setminus\{\0\}$.

Let $C\subset \R^n$ be an irreducible variety of dimension $d$, $1\le d <n$ as in \S\ref{sec:irredvar}.  In this section we show that, as in the case of $\ell_2$ norm, $\delta_{\|\cdot\|_p} (C)$ 
in Definition \ref{defnnudegC}  is the degree of $\pi: V_{1,\C}\to \C^n$, and $V_{1,\C}$  can be considered as the   critical point correspondence.
(That is, in Definition \ref{defnnudegC} we take $k=1$.)

Let $g(\x):=\|\x\|_p^p$. 
Assume that $\y\in C_\C\setminus\Sing C_\C$.  Denote by $DC(\y)\in\C^{n\times m}$ the complex Jacobian as defined in \S\ref{sec:irredvar}.  Let $\U(\y)\subset \C^n$ be 
the subspace spanned by the columns  of $DC(\y)$.  Then $\dim \U(\y)=n-d$.  
Let 
\[\Sigma_{2m,2l+1,1}(C)=\{(\z,\y)\in \C^n\times (C_\C\setminus \Sing C_\C), \quad \bF_{2m-2l-1}(\z-\y)\in \bF_{2l+1}(\U(\y))\}.\]

That  is, 
$$\Sigma_{2m, 2l+1,1}(C)=\{(\z,\y)\in \C^n\times (C_\C\setminus \Sing C_\C), \  \z\in\y+\bF_{2m-2l-1}^{-1}(\bF_{2l+1}\U(\y))\}.$$\\
Let further 
$\Sigma_{2m, 2l+1}(C)=\textrm{Closure} (\Sigma_{2m, 2l+1,1}(C))$.
Lemma \ref{charcritpoits} yields that any critical point of $\y_0$ of  $g_\x(\y):=\|\x-\y\|_p^p$ in $C\setminus \Sing C$ satisfies $(\x,\y_0)\in \Sigma_{2m,2l+1,1}(C)$. 
\begin{lemma}\label{sig2malgvar}  Let $m>l\ge 0$ be integers.  Then $\Sigma_{2m, 2l+1}(C)\subset \C^n\times \C^n$ is a closed algebraic variety, 
and each of its irreducible components is of dimension $n$.
Furthermore, there exists exactly one irreducible component $\Sigma_{2m,2l+1}^1(C)$ of $\Sigma_{2m, 2l+1}(C)$ containing all real points in $\Sigma_{2m, 2l+1}(C)$.
\end{lemma}
\proof  Consider first the case when $m=1, l=0$.   Then $\Sigma_{2,1}(C)=\Sigma_{\|\cdot\|^2}(C)_\C$ is the critical correspondence variety studied in \cite{DHOST,FS13},
and discussed in \S\ref{sec:irredvar}.   
Observe first that 
\[\Phi(C):=\{(\uu,\y),\; \y\in C_\C\setminus \Sing C_\C, \uu\in \U(\y)\}\] 
is a quasi-algebraic variety, (see Lemma \ref{charcritpoits}).  Furthermore, it 
is isomorphic to an $(n-d)$-dimensional vector bundle over $C_\C\setminus\Sing C_\C$.   Hence it is a connected complex manifold of dimension $n$.
Therefore Closure$(\Phi(C))$ is an irreducible complex variety of dimension $n$. (It is usually called the Nash modification of $C_{\C}$.)
Consider the following linear automorphism $A$ of $\C^{n}\times\C^{ n}$:
$(\z,\w)\mapsto (\z+\w,\w)$.  Then 
\begin{eqnarray*}
&&A(\Sigma_{2,1,1}(C))=\Phi(C)\\
&&\tilde \Sigma_{2,1}(C):=\textrm{Closure}(\Phi(C))= A(\Sigma_{2,1}(C)).
\end{eqnarray*}
Hence $\tilde \Sigma_{2,1}(C)$ and $\Sigma_{2,1}(C)$  are $n$-dimensional irreducible varieties.

Let $\tilde \bF_{2m-1}:\C^n\times \C^n$ be given by $(\z,\w)\mapsto (\bF_{2m-1}(\z),\w)$.  Clearly, $\tilde\bF_{2m-1}$ is a proper polynomial map
of degree $(2m-1)^n$.  Hence, for each irreducible variety in $W\subset \C^n\times\C^n$,   $\tilde\bF_{2m-1}^{-1}(W)$ is a union of at most $(2m-1)^n$ irreducible
varieties, each of dimension $\dim W$.  Furthermore, $\tilde \bF_{2m-1}(W)$ is an irreducible variety of dimension $n$.  Let
\begin{equation}\label{defSig2m'}
\tilde\Sigma_{2m,2l+1}(C)=\tilde \bF_{2m-2l-1}^{-1}(\tilde\bF_{2l+1}(\tilde\Sigma_{2,1}(C)))=\bigcup_{i=1}^{N(2m)} \tilde\Sigma_{2m,2l+1}^{i}(C).
\end{equation}
Here  $\tilde\Sigma_{2m,2l+1}^{i}(C)$ are distinct irreducible components of $\tilde\Sigma_{2m,2l+1}(C)$ for $i=1,\ldots,N(2m)$.
It is straightforward to show
\begin{eqnarray}\notag
&&\Sigma_{2m,2l+1}(C)=A^{-1}(\tilde\Sigma_{2m,2l+1}(C))=\bigcup_{i=1}^{N(2m)} \Sigma_{2m,2l+1}^i(C)\\  
&&\Sigma_{2m,2l+1}^i(C):=A^{-1}(\tilde\Sigma_{2m,2l+1}^i(C)), \;i=1,\ldots,N(2m).\label{Sigma2mchar}
\end{eqnarray}
Here  $\Sigma_{2m,2l+1}^i(C)$ are distinct irreducible components of $\Sigma_{2m,2l+1}(C)$ of dimension $n$ for $i=1,\ldots,N(2m)$.

Assume now that $(\x,\y)$ is a real point in $\tilde\Sigma_{2m,2l+1}(C)$.  So $\uu:=\bF_{2m-2l-1}(\x)\in\bF_{2l+1}(\U(\y))\cap\R^n$. 
 Clearly, $\bF_{2m-2l-1}^{-1}(\uu)$ contains exactly one
real point $\x$.  Hence $\tilde\Sigma_{2m,2l+1}(C)\cap (\R^n\times C)=\tilde\Sigma_{2m,2l+1}^1(C) \cap (\R^n\times C)$.   
Thus $\Sigma_{2m,2l+1}^1(C)$ of $\Sigma_{2m,2l+1}(C)$ contains all real points in $\Sigma_{2m,2l+1}(C)$.  \qed

\begin{theo}\label{finitnubcritpts}  Let $m>l\ge 0$ be integers and assume that $\Sigma_{2m,2l+1}(C)$ is defined as above.
Let $\Sigma_{2m,2l+1}^1(C)$ be the irreducible component of dimension $n$ of $\Sigma_{2m,2l+1}(C)$ which contains all the real points of  $\Sigma_{2m,2l+1}(C)$.
Let $p=\frac{2m}{2l+1}$ and  $g(\x)=\|\x\|_p^p$.  Let $\Omega_{\|\cdot\|_p}(C)$ and $\Sigma_{g,1}(C)$ be defined as in Corollary   
\ref{defOmegnuC} and Lemma \ref{theo:CritSemialg} respectively.  Then $\Omega_{\|\cdot\|_p}(C)\cap \Sigma_{g,1}(C)\subset \Sigma_{2m,2l+1}^1(C)$.
Hence $\delta_{\|\cdot\|_p}(C)$ is the degree of the dominating map $\pi:\Sigma_{2m,2l+1}^1(C)\to\C^n$.  That is, $\Sigma_{2m,2l+1}^1(C)$ is the
critical point correspondence for  the $\ell_p$ norm.

\end{theo}
\proof  As in the proof of Lemma \ref{interseclem}, we have 
\[\Omega_{\|\cdot\|_p}(C)\setminus(\R^n\times \Sing C)=\Omega_{\|\cdot\|_p}(C)\cap \Sigma_{g,1}(C)
\subset \Sigma_{2m,2l+1,1}(C)\subset \Sigma_{2m,2l+1}(C).\]  
As  $\Omega_{\|\cdot\|_p}(C)\setminus (\R^n\times \Sing C)\subset\R^n$
it follows that $\Omega_{\|\cdot\|_p}(C)\setminus(\R^n\times \Sing C)\subset \Sigma_{2m,2l+1}^1(C)$. 
Let  $\Omega_\nu(C)\setminus(\R^n\times \Sing C)=\bigcup_{i=1}^N \Phi_i$ be a smooth  Whitney semi-algebraic stratification as in Theorem \ref{theodefdegnuC}. 
Then each $\Phi_i$ of dimension $n$ is an open semi-algebraic subset of $\Sigma_{2m,2l+1}^1(C)$.  Theorem \ref{theodefdegnuC} yields that 
 $\pi:\Sigma_{2m,2l+1}^1(C)\to\C^n$ is dominating.   Definition \ref{defnnudegC} yields that  $\delta_{\|\cdot\|_p}(C)$ is the degree of the dominating map $\pi:\Sigma_{2m,2l+1}^1(C)\to\C^n$.\qed

 \section{Algebraic critical point correspondences}\label{sec:algcritptcor}

In this section we give more detailed information about the variety $V_{i,\C}$ appearing in Theorem \ref{theodefdegnuC}.
Let $\nu$ be a semi-algebraic norm such that $\nu$ and $\nu^*$ are differentiable.
Consider the smooth Whitney stratification of $\nu:\R^n\to \R$ given by Proposition \ref{prop:strata}.
Assume that $\Delta_i$ is of dimension $n$.  (Since $\nu$ is not differentiable at $\0$, we deduce that $\0\not\in \Delta_i$.)
Then the graph of $\nu \mid \Delta_i$ is an open semi-algebraic set in an irreducible variety $W_i\subset \R^n\times \R$
of dimension $n$.  Let $W_{i,\C}\subset \C^n\times \C$ be the complexification of $W_i$.   Hence $W_{i,\C}$ is a hypersurface in $\C^{n}\times \C$.
For simplicity of notation we let $W=W_{i,\C}$ when no ambiguity would arise.   Thus 
\[W=\{(\z,t)\in\C^{n+1}, \;G(\z,t)=0\},\]
 where
\begin{equation}\label{defGzt}
G(\z,t)=\sum_{j=0}^M a_j(\z)t^j, \quad a_j(\z)\in \C[\z], \; j=0,\ldots,M.
\end{equation}
The arguments below show that $G(\z,t)$ is an irreducible polynomial in $\C[(\z,t)]$.

Let $\tau:\C^n\times \C\to \C^n$ be the projection onto the first component.  Since $\tau(W_i)\supset \Delta_i$, it follows that $\tau: W_{i,\C}\to \C^n$ is a  
dominating polynomial map of degree $m$. Furthermore,  
 $W_{i,\C}$ is the graph of a multivalued algebraic function $f_i$.  For simplicity of notation we also let $f:=f_i$ when no ambiguity would arise.
That is,  $f$ satisfies a polynomial equation $G(\z,f)=0$ where $G(\z,t)$ is given by \ref{defGzt}.
Let 
\[W':=\{(\z,t)\in W,\;\frac{\partial G}{\partial t}=0\}.\] 
Clearly, $W'$ is a strict subvariety of $W$.  For $(\z_0,t_0)\in W\setminus W'$ the implicit function theorem yields that  $f$ is an analytic function of $\z$ 
  in a neighborhood of $\z_0$ with $f(\z_0)=t_0$.
Let  
\[\bF:=(\frac{\partial G}{\partial z_1},\ldots,\frac{\partial G}{\partial z_n}).\]
Note that $\bF:\C^{n+1}\to \C^n$ is a polynomial map $(\z,t)\mapsto \bF((\z,t))$ and  that 
\begin{equation}\label{nuzderfor}
\nabla f(\z)=\bH((\z,t)):=-\frac{\bF((\z,t))}{\frac{\partial G(\z,t)}{\partial t}}, \quad (\z,t)\in W\setminus W'.
\end{equation}

For $\y\in C_\C\setminus\Sing C_\C$ let $\U(\y)$ be defined as in Lemma \ref{charcritpoits}. (Recall that $C$ is a $d$-dimensional variety in $\R^n$.)
The critical set $\Sigma_{\nu,1}(C,W)$ is defined as 
\[
\Sigma_{\nu,1}(C,W)=\{(\z,t,\y)\in \C^n\times \C\times (C_\C\setminus \Sing C_\C), (\z-\y,t)\in W,\bF(\z-\y,t)\in \U(\y)\},\\
\]
The arguments for equivalence of the conditions in Lemma  \ref{charcritpoits} yield that $\Sigma_{\nu,1}(C,W)$ is a quasi-variety.  
Hence $\Sigma_\nu(C,W)=\textrm{Closure}(\Sigma_{\nu,1}(C,W))\label{defSigmagCW}$ is an algebraic set.
As in the case of the $\ell_{2m}$ norm,  $\Sigma_\nu(C,W)$ may contain  several irreducible components.
\begin{theo}\label{factorSigma}  
Assume that $\nu,\nu^*\in C^1(\R^n\setminus\{\0\})$ and $\nu$ is semi-algebraic.
Consider the  smooth Whitney stratification of $\nu:\R^n\to \R$ given in Proposition \ref{prop:strata}.
Assume that $\Delta_i$ is of dimension $n$.  Let $W_i\subset \R^n\times \R$ be an $n$-dimensional irreducible variety that contains the graph of $\nu\mid \Delta_i$.
Let  $W_{i,\C}$ denote the complexification of $W_i$.
Let $C\subset \R^n$ be an irreducible variety of dimension $d$.
With $\Sigma_\nu(C,W_{i,\C})$ defined as above and and $\Omega_\nu(C)$ as in  Corollary \ref{defOmegnuC}, respectively, define  

\[\Omega_\nu(C,\Delta_i):=\{(\x,\z)\in\Omega_\nu(C)\setminus (\R^n\times \Sing C), \; \x-\z\in \textrm{ Closure}(\Delta_i)\setminus\{\0\}\}.\]

Then $\Omega_\nu(C,\Delta_i)$ is semi-algebraic. Furthermore, 
\begin{enumerate}
\item Assume that $\pi(\Omega_\nu(C,\Delta_i))$ has dimension $n$.
Then $\Sigma_\nu(C,W_{i,\C})$ contains a positive number of subvarieties $\tilde V_{i,1},\ldots,\tilde V_{i,l_i}$
such that $\dim \tilde V_{i,j}=n$  for $j=1,\ldots,l_i$.
\item Let $\omega:\C^n\times \C\times \C^n \to \C^n\times \C^n$ be the projection onto the product of the first and the last factor. 
Then for each $j\in \{1,\ldots,l_i\}$ there exists a variety $V_{j'}$ in the set of $k$ varieties $\{V_1,\ldots,V_k\}$, as given in Definition \ref{defnnudegC},
such that $\omega: \tilde V_{i,j}\to V_{j'}$ is dominating. 
\end{enumerate}
Moreover, for each $V_l\in\{V_1,\ldots,V_k\}$, there exists $i$ and $j$ such that $\omega: \tilde V_{i,j}\to V_l$.
\end{theo}
\proof 
Since $\Omega_\nu(C)$ and $\Delta_i$ are semi-algebraic, we deduce immediately that $\Omega_\nu(C,\Delta_i)$ is semi-algebraic.
Assume that the semi-algebraic set  $\Delta'_i:=\pi(\Omega_\nu(C,\Delta_i))$ has dimension $n$.
Denote by $G_i(\z,t)$ the polynomial induced by $\nu\mid \Delta_i$.
Let
\[\tilde \Omega_\nu(C,\Delta_i):=\{(\x,t,\z),  (\x,\z)\in \Omega_\nu(C,\Delta_i), t=\nu(\x-\z)\}.\]
We claim that $\tilde \Omega_\nu(C,\Delta_i)\subset \Sigma_\nu(C,W_{i,\C})$.  
Since $W_i$ is a closed set, it contains the graph of $\nu\mid \textrm{Closure}(\Delta_i)$.
Assume that $(\x,t,\z)\in \tilde\Omega_\nu(C,\Delta_i)$.  By the definition, $(\x-\z,t)\in W_{i,\C}$.  Furthermore, $\z\in C$ is a unique best $\nu$-approximation to $\x$.
By definition, $\x-\z\in \textrm{ Closure}(\Delta_i)\setminus \{\0\}$.  Hence $\x-\z\ne \0$.  Therefore $\nabla\nu(\x-\z)\ne\0$.
As $(\x-\z,\nu(\x-\z))\in W_{i}$ it follows that $\bF(\x-\z,\nu(\x-\z))\in \U(\z)$.  Thus $(\x,\nu(\x-\z),\z)\in \Sigma_\nu(C,W_{i,\C})$.

Let  $\Omega_\nu(C)\setminus (\R^n\times \Sing C)=\bigcup_{l=1}^N \Theta_l$ be the  smooth Whitney  decomposition as in Theorem \ref{theodefdegnuC}.
Assume that $\dim \Theta_l=n$.  Hence $\dim \pi(\Theta_l)=n$.  Suppose furthermore that $\dim \pi(\Theta_l)\cap \Delta'_i=n$.  Then the variety 
$Y_{i}=$Closure $(\omega(W_{i,\C}))$ contains
$\Omega_\nu(C,\Delta_i)$.   Let $V_{l,\C}$ be the irreducible variety defined in Theorem \ref{theodefdegnuC}.  Since $V_{l}$ is the minimal variety 
containing $\Theta_l$, it follows that $V_{l,\C}\subset Y_i$.   Let
\[\tilde Y_{i,l}:=\{(\x,t,\z)\in \Sigma_\nu(C,W_{i,\C}), \;(\x,\z)\in V_{l,\C}\}.\]
Let $\tilde V_{i,l}\subset \tilde Y_{i,l}$ be the smallest subvariety of $\tilde Y_{i,l}$ which contains the points of the semi-algebraic set $\tilde \Omega_\nu(C,\Delta_i)$
of dimension $n$.  Then $\dim \tilde V_{i,l}=n$ and $\omega: \tilde V_{i,l}\to V_{l,\C}$ is dominating.

It remains to show that each $V_l$ appearing in Theorem \ref{theodefdegnuC} corresponds to some $\Theta_l$ of dimension $n$.
Clearly, $\diag(C):=\{(\z,\z), \z\in C\}$ is an irreducible variety of dimension $d<n$.  By refining the stratification if needed we can assume that 
in the decomposition of $\Omega_\nu(C)\setminus (\R^n\times \Sing C)$ each $\Theta_i$ of dimension $n$ does not intersect $\diag(C)$,
i.e. $\Theta_i\cap \diag(C)=\emptyset$ for $i=1,\ldots,M$.

Assume that in the smooth Whitney stratification of $\nu:\R^n\to \R$ given in Proposition \ref{prop:strata} we have that $\dim \Delta_i=n$ for $i=1,\ldots,L'$
and $\dim \Delta_i<n$ for $i>L'$.  Hence  $\R^n\setminus (\bigcup_{i=1}^{L'} \Delta_i)$ is a semi-algebraic set of dimension  at most $n-1$.  Therefore
\[\R^n=\bigcup_{i=1}^M \textrm{Closure} (\Delta_i), \quad \R^n\setminus\{\0\}=\bigcup_{i=1}^M \textrm{Closure} (\Delta_i)\setminus\{\0\}.\]

Assume that $\dim\Theta_l=n$.  So $\dim\pi(\Theta_l)=n$.  As $\Theta_l\cap \diag(C)=\emptyset$ it follows that $\x-\z\ne 0$ for each $(\x,\z)\in \Theta_l$.
Hence 
\[\Theta_l=\cup_{i=1}^{L'} (\Theta_l\cap\Omega_\nu(C,\Delta_i)), \quad \pi(\Theta_l)=\cup_{i=1}^{L'} (\pi(\Theta_l\cap\Omega_\nu(C,\Delta_i))).\]
Therefore there exists $\Delta_i$ of dimension $n$ such that $\dim \pi(\Theta_l)\cap\pi(\Omega_\nu(C,\Delta_i))=n$.\qed

\end{document}